\documentclass[12pt]{article}
\usepackage{cite}
\usepackage{amssymb,amsmath,latexsym}
\usepackage{amsthm}

\usepackage{pdfsync}
\usepackage[matrix,arrow,curve]{xy}
\input xy
\xyoption{all}

\oddsidemargin=\evensidemargin \advance\topmargin by-\headheight

\def\VRT#1{*=<7mm>[o][F-]{#1}}

\DeclareMathOperator{\aut}{Aut}

\DeclareMathOperator{\cyc}{Cyc}

\DeclareMathOperator{\hol}{Hol}

\DeclareMathOperator{\orb}{Orb}

\DeclareMathOperator{\pr}{pr}
\DeclareMathOperator{\rad}{rad}
\DeclareMathOperator{\rk}{rk}

\DeclareMathOperator{\Span}{span}

\DeclareMathOperator{\sym}{Sym}

\def\mZ{{\mathbb Z}}

\def\A{{\cal A}}

\def\H{{\cal G}}

\def\M{{\cal M}}

\def\S{{\cal S}}

\def\twoe{\underset{\scriptscriptstyle ^2}{\approx}}

\def\lg{\langle}

\def\rg{\rangle}

\def\proof{\noindent {\bf Proof}.\ }
\def\bull{\vrule height .9ex width .8ex depth -.1ex }

\def\qaq{\quad\text{and}\quad}

\makeatletter
\renewcommand{\subsection}{\@startsection{subsection}{2}{0mm}{-2mm}{-2mm}
{\bf\normalsize}}
\def\sbsn{\subsection{\hspace{-3mm}}}

\makeatother

\newtheorem{formula}{}[section]
\newtheorem{proposition}[formula]{Proposition}
\newtheorem{definition}[formula]{Definition}
\newtheorem{corollary}[formula]{Corollary}
\newtheorem{lemma}[formula]{Lemma}
\newtheorem{theorem}[formula]{Theorem}

\newtheorem{hypothesis}[formula]{Conjecture}
\newtheorem{example}[formula]{Example}
\newtheorem{question}[formula]{Question}

\newtheorem*{thrmP}{Theorem (P\"oschel)}
\theoremstyle{remark}
\newtheorem{remark}[formula]{Remark}

\def\thrm{\begin{theorem}}
\def\thrml#1{\begin{theorem}\label{#1}}
\def\ethrm{\end{theorem}}
\def\prpstn{\begin{proposition}}
\def\prpstnl#1{\begin{proposition}\label{#1}}
\def\eprpstn{\end{proposition}}
\def\rmrk{\begin{remark}}
\def\rmrkl#1{\begin{remark}\label{#1}}
\def\ermrk{\end{remark}}
\def\dfntn{\begin{definition}}
\def\dfntnl#1{\begin{definition}\label{#1}}
\def\edfntn{\end{definition}}
\def\nmrt{\begin{enumerate}}
\def\enmrt{\end{enumerate}}
\def\tm#1{\item[{\rm (#1)}]}
\def\qtn{\begin{equation}}
\def\qtnl#1{\begin{equation}\label{#1}}
\def\eqtn{\end{equation}}
\def\lmm{\begin{lemma}}
\def\lmml#1{\begin{lemma}\label{#1}}
\def\elmm{\end{lemma}}
\def\crllr{\begin{corollary}}
\def\crllrl#1{\begin{corollary}\label{#1}}
\def\ecrllr{\end{corollary}}
\def\hpthss{\begin{hypothesis}}
\def\hpthssl#1{\begin{hypothesis}\label{#1}}
\def\ehpthss{\end{hypothesis}}
\def\xmpl{\begin{example}}
\def\xmpll#1{\begin{example}\label{#1}}
\def\exmpl{\end{example}}
\def\css{\begin{cases}}
\def\ecss{\end{cases}}
\def\qstn{\begin{question}}
\def\qstnl#1{\begin{question}\label{#1}}
\def\eqstn{\end{question}}

\begin{document}

\title{On schurity of finite abelian groups}
\author{
Sergei Evdokimov \\[-2pt]
\small St.Petersburg Department of Steklov Mathematical Institute\\[-4pt]
\small Fontanka 27, St. Petersburg, Russia \\[-4pt]
{\tt \small evdokim@pdmi.ras.ru }
\thanks{The work was partially supported by the Slovenian-Russian bilateral project, grant no.~BI-RU/12-13-035.
The second author was also supported by the project ARRS -
Agencija za raziskovanje Republike Slovenija, program no.~P1-0285. The third author  was partially supported by the RFFI Grant 11-01-00760-a.}
\and
Istv\'an Kov\'acs\\[-2pt]
\small IAM and FAMNIT, University of Primorska \\ [-4pt]
\small Muzejski trg 2, 6000 Koper, Slovenia\\[-4pt]
{\tt \small istvan.kovacs@upr.si}
\and
Ilya Ponomarenko\\[-2pt]
\small St.Petersburg Department of Steklov Mathematical Institute\\ [-4pt]
\small Fontanka 27,  St. Petersburg, Russia \\[-4pt]
{\tt \small inp@pdmi.ras.ru}
}

\maketitle

\begin{abstract}
A finite group $G$ is called a Schur group, if any Schur ring over~$G$ is
associated in a natural way with a subgroup of $\sym(G)$ that contains all right
translations. Recently, the authors have completely identified the cyclic Schur groups.
In this paper it is shown that any abelian Schur group belongs to one of several
explicitly given families only.
In particular, any non-cyclic abelian Schur group
of odd order is isomorphic to $\mZ_3\times\mZ_{3^k}$ or $\mZ_3\times\mZ_3\times\mZ_p$ where $k\ge1$ and
$p$ is a prime.  In addition, 
we prove  that
$\mZ_2\times\mZ_2\times\mZ_p$ is a Schur group for every prime~$p$.  \medskip
\end{abstract}

\section{Introduction}

Let $G$ be a finite group.   A subring of the group ring $\mZ G$ is called a {\em Schur ring} ({\em S-ring} for short) over $G$ if it
is determined in a natural way by a special partition of $G$ (for the exact definition we refer to
Section~\ref{120813a}). The definition goes back to Schur and Wielandt, they
used "the S-ring method"\ to study a permutation group having a regular subgroup, see~\cite{Sch33,Wie64}. Namely,
with any group $\Gamma\le\sym(G)$ containing all right (or left) translations, one can associate the S-ring determined by the
partition of $G$ into the orbits of the stabilizer of the point $1_G$ in $\Gamma$. Following~\cite{Poe74} any S-ring obtained in this way is called
{\it schurian}; an explanation of the term is that Schur himself studied only such S-rings. \medskip

Wielandt wrote in~\cite[p.54]{Wie69} that ``Schur had conjectured for a long time that every S-ring is determined by a suitable permutation group'',
in our terms  that every S-ring is schurian. However, this proved to be false and first counterexamples were found
in~\cite[Theorem~25.7]{Wie64}: these examples form a family of non-schurian S-rings over the  group $G=\mZ_p\times\mZ_p$
where $p \ge 5$  is a prime. In terms of the following definition given by P\"oschel in~\cite{Poe74},
this result shows that $G$ is not a Schur group.

\dfntn
A finite group $G$ is called a Schur group if  every S-ring over $G$ is schurian
\edfntn

P\"oschel proved in the same paper that any section of a Schur group is also Schur, and asked for
the structure  of Schur groups. He was able himself to characterize Schur $p$-groups as follows.

\begin{thrmP}
Given a prime $p\ge 5$ a finite $p$-group is Schur if and only if it is cyclic.
\end{thrmP}

The "if"\ part of the theorem states that any cyclic $p$-group is Schur for $p\ge 5$. In fact, in \cite{Poe74} it
was also proved that any cyclic $3$-group is Schur. Later the second author proved the same statement for $p=2$ in~\cite{Kov11}.
Taking into account that Schur himself mainly worked with S-rings over a cyclic group and that a cyclic group
the order of which is  the product of two distinct primes is also Schur~\cite{KliP81}, it was not suprising that 
some experts beleived  
that any  finite cyclic group is Schur (the Schur-Klin
conjecture). However, this conjecture was disproved in  \cite{EvdP01}. Then further examples of
cyclic Schur and non-Schur groups were found  in \cite{EvdP12}. Eventually, the following characterization
of cyclic Schur groups was recently completed by the  authors  in \cite{EvdKP}.\medskip

{\bf \noindent Cyclic Schur Group Theorem.}\ {\it
A cyclic group of order $n$ is a Schur group if and only if $n$ belongs to one of the following
families of integers:
$$ p^k, \; pq^k, \; 2pq^k, \; pqr, \; 2pqr,$$
where $p,q,r$ are distinct primes, and $k \ge 0$ is an integer. }\medskip

In the present paper we are interested in the P\"oschel question for abelian groups.  First we prove a sufficient condition
for an abelian group to be non-Schur (Theorem~\ref{040211b}). The proof is based on the construction of a non-schurian
S-ring as the generalized wreath product of normal cyclotomic S-rings; such a product was studied in \cite{EvdP01,EvdP12}
and used in~\cite{EvdKP} to prove the Cyclic Schur Group Theorem.
This condition is used to prove Theorems~\ref{290713a} and~\ref{040213a}  showing  that any non-cyclic abelian Schur group belongs
to a relatively short list of groups. Here we also apply the enumeration of  small\footnote{Over a group of order at most $41$.}
S-rings 
that was done by  C. Pech and S. Reichard
(the computation results were announced in \cite[Section~8.4]{PecR09}). The obtained catalog is used
throughout the paper, and will be refered  to as the  {\em Classification of Small S-rings} (CSS).

\thrml{290713a}
An elementary abelian non-cyclic group of order $n$ is Schur if and only if $n\in\{4,8,9,16,27,32\}$.
\ethrm

\thrml{040213a}
An abelian Schur group which is neither cyclic  nor elementary abelian,  is isomorphic to one of the following
nine families of groups:
\nmrt
\tm{1} $\mZ_2 \times \mZ_{2^k},$ $\mZ_{2p}\times\mZ_{2^k}$, $E_4\times\mZ_{p^k}$, $E_4\times\mZ_{pq}$, $E_{16}\times\mZ_p$,
\tm{2} $\mZ_3 \times \mZ_{3^k},$ $\mZ_6\times\mZ_{3^k}$,  $E_9\times\mZ_q,$ $E_9\times\mZ_{2q}$.
\enmrt
where $p$ and $q$ are distinct primes, $p\ne 2$, and $k \ge 1$ is an integer.
\ethrm

The following statement immediately follows from Theorems~\ref{290713a} and~\ref{040213a}.

\crllr
Any non-cyclic abelian Schur group of odd order is isomorphic to either $\mZ_3\times\mZ_{3^k}$ or $E_9\times\mZ_p$ where $p$ is a prime. 
\ecrllr

We do not know  whether  Theorem~\ref{040213a} can be conversed. However, by  the CSS  this is true for any group of order at most~$41$.
Moreover, the following result shows that this is also true for groups belonging to the intersection of the second and third families.

\thrml{040213c}
For any prime $p$ the group $E_4 \times \mZ_p$ is a Schur one.
\ethrm

The above discussion naturally leads to the following question.

\qstnl{040213b}
Is it true that all groups cited  in  Theorem~\ref{040213a} are Schur?
\eqstn

To make the paper self-contained we collect the basic facts on S-rings  in Section~\ref{120813a}. In
Section~\ref{140813a} we discuss the concepts concerning schurian S-rings.\medskip

{\bf Notation.}
As usual by $\mZ$ we denote the ring of rational integers.

The set of non-identity elements of a group $G$ is denoted by  $G^\#$.

For a normal subgroup $H$ of $G$ the quotient epimorphism from $G$ onto $G/H$ is denoted by $\pi_{G/H}$.

The subgroup of $G$ generated by a set $X\subset G$ is denoted by $\lg X\rg$; we also set
$
\rad(X)=\{g\in G:\ gX=Xg=X\}.
$

The group of all permutations of $G$ is denoted by $\sym(G)$.

For a set $\Delta\subset\sym(G)$ and a section $S$ of  $G$ we set
$$
\Delta^S=\{f^S:\ f\in \Delta,\ S^f=S\},
$$
where $f^S$ is the bijection of $S$ induced by~$f$.

The subgroup of $\sym(G)$ induced by right multiplications of ~$G$  is denoted by $G_{right}$.

The holomorph $\hol(G)$ is identified with the subgroup of~$\sym(G)$ generated by $G_{right}$ and~$\aut(G)$.

The orbit set of a group $\Gamma\le\sym(G)$ is denoted by $\orb(\Gamma)=\orb(\Gamma,G)$.

We write $\Gamma\twoe\Gamma'$ if groups $\Gamma,\Gamma'\le\sym(G)$ are $2$-equivalent, i.e.
have the same orbits in the coordinate-wise action on~$G\times G$.

The cyclic group of order $n$ is denoted by  $\mZ_n$.

The elementary abelian group of order $p^k$ is denoted by $E_{p^k}$.

\section{S-rings}\label{120813a}
In what follows regarding S-ring theory we use the  notations and terminology from~\cite{EvdP12}.\medskip

Let $G$ be a finite group. A subring~$\A$ of the group ring~$\mZ G$ is called a {\it Schur ring} ({\it S-ring}, for short) over~$G$ if 
there exists a partition $\S=\S(\A)$ of~$G$ such that
\nmrt
\tm{S1} $\{1_G\}\in\S$,
\tm{S2} $X\in\S\ \Rightarrow\ X^{-1}\in\S$,
\tm{S3} $\A=\Span\{\sum_{x\in X}x:\ X\in\S\}$.
\enmrt
A group isomorphism $f:G\to G'$ is called a {\it Cayley isomorphism} from an S-ring $\A$ over $G$ to an S-ring $\A'$ over $G'$ if $\S(\A)^f=\S(\A')$.\medskip

The elements of $\S$ and the number $\rk(\A)=|\S|$ are called respectively the {\it basic sets} and the {\it rank} of the S-ring~$\A$. Any union of basic sets
is called an {\it $\A$-subset of~$G$} or {\it $\A$-set}; the set of all of them is denoted by $\S^*(\A)$. The
latter set is closed with respect to taking inverse and product. Given $X\in\S^*(\A)$ the submodule of~$\A$ spanned by the set
$$
\S(\A)_X=\{Y\in\S(\A):\ Y\subset X\}
$$
is denoted by $\A_X$.\medskip  

Any subgroup of $G$ that is an $\A$-set is called an {\it $\A$-subgroup} of~$G$ or {\it $\A$-group}; the set of all of them is denoted by $\H(\A)$.
With each $\A$-set $X$ one can  naturally associate two $\A$-groups, namely $\lg X\rg$ and  $\rad(X)$ (see Notation).
The following useful lemma was proved in \cite[p.21]{EvdP09}.

\lmml{090608a}
Let $\A$ be an S-ring over a group $G$, $H\in\H(\A)$ and $X\in\S(\A)$. Then
the cardinality of the set $X\cap Hx$ does not depend on $x\in X$.
\elmm

Let $S=U/L$ be a  section  of $G$. It is  called an  {\it $\A$-section}, if both  $U$ and $L$ are $\A$-groups. Given $X\in\S(\A)_U$
 the module
$$
\A_S=\Span \{\pi_S(X):\ X\in\S(\A)_U\}
$$
is an S-ring over the group~$S$, the basic sets of which are exactly the sets from the right-hand side of the formula.

\dfntn
The S-ring $\A$ is called an {\it $S$-wreath product} if the group $L$ is
normal in~$G$ and $L\le\rad(X)$ for all basic sets $X$ outside $U$; in this case we write
\qtnl{050813a}
\A=\A_U\wr_S\A_{G/L},
\eqtn
and omit $S$ when $U=L$.
\edfntn

When the explicit indication of the section~$S$ is not important, we use the term {\it generalized wreath product}.
The $S$-wreath product is {\it nontrivial} or {\it proper}  if $1\ne L$ and $U\ne G$.\medskip

It was proved in \cite[Theorem~3.1]{EvdP01}  that given
a section $S=U/L$ of an abelian group~$G$, and S-rings $\A_1$ and $\A_2$ over the groups $U$ and $G/L$ respectively
such that $S$ is both an $\A_1$- and an $\A_2$-section, and
$$
(\A_1)_S=(\A_2)_S,
$$
there is a uniquely determined $S$-wreath product~\eqref{050813a} such that $\A_U=\A_1$ and $\A_{G/L}=\A_2$.\medskip

If $\A_1$ and $\A_2$ are S-rings over groups $G_1$ and $G_2$ respectively, then the subring $\A=\A_1\otimes \A_2$ of the ring $\mZ G_1\otimes\mZ G_2=\mZ G$
where $G=G_1\times G_2$, is an S-ring over the group $G$ with
$$
\S(\A)=\{X_1\times X_2: X_1\in\S(\A_1),\ X_2\in\S(\A_2)\}.
$$
It is called the {\it tensor product} of $\A_1$ and $\A_2$.

\lmml{050813b}
Let $\A$ be an S-ring over an abelian group $G=G_1\times G_2$. Suppose that $G_1$ and $G_2$ are $\A$-groups. Then
\nmrt
\tm{1} $\pr_{G_i}(X)\in\S(\A)$ for all $X\in\S(\A)$ and $i=1,2$,
\tm{2} $\A\ge\A_{G_1}\otimes \A_{G_2}$,  and the equality is attained whenever $\A_{G_i}=\mZ G_i$ for some $i\in\{1,2\}$.
\enmrt
\elmm
\proof Statement~(1) and the inclusion in statement~(2) were proved in \cite[Lemma~2.2]{EvdP09}. To prove the rest 
without loss of generality we can assume that $\A_{G_1}=\mZ G_1$ and consequenly $\S(\A_{G_1})$ consists of singletons.
Then obviously $\{x_1\}\times X_2\in\S(\A)$ for all $x_1\in G_1$ and $X_2\in\S(\A_{G_2})$.
Since the union of the above products forms a partition of $G$, we are done.\bull\medskip

The following important theorems go back to Schur and Wielandt (see \cite[Ch.~IV]{Wie64}). The first of them is known as the
Schur theorem on multipliers, see~\cite{EvdP09}.

\thrml{261009b}
Let $\A$ be an S-ring over an abelian group $G$. Then any central element of $\aut(G)$ is a Cayley automorphism of~$\A$.
\ethrm

In general, Theorem~\ref{261009b} is not true when an automorphism of $G$ is replaced by a homomorphism of $G$ to itself.
However, the following weaker statement holds.

\thrml{261009w}
Let $\A$ be an S-ring over an abelian group $G$. Then given a prime $p$ dividing $|G|$, the set
\qtnl{030713a}
X^{[p]}:=\{x^p:\ x\in X,\ |xH\cap X|\not\equiv 0\pmod p\}
\eqtn
where $H=\{g\in G:\ g^p=1\}$, belongs to $\S^*(\A)$ for all $X\in\S(\A)$.
\ethrm

\section{Schurian S-rings and automorphism groups}\label{140813a}
Let $G$ be a finite group. It was proved by Schur (see \cite[Theorem~24.1]{Wie64}) that  any group
$\Gamma\le\sym(G)$ that contains $G_{right}$ produces an S-ring~$\A$ over $G$ such that
$$
\S(\A)=\orb(\Gamma_1,G)
$$
where $\Gamma_1=\{\gamma\in\Gamma:\ 1^\gamma=1\}$ is the stabilizer of the point $1=1_G$ in $\Gamma$. Any such  S-ring is called {\it schurian}. 
Group rings and S-rings of rank~$2$ are obviously schurian. 
A less trivial example is given by a {\it cyclotomic} S-ring $\A=\cyc(K,G)$ where $K\le\aut(G)$; by definition  $\S(\A)=\orb(K,G)$.\medskip

The schurity concept is closely related to  automorphisms of an S-ring.
Nowadays there are several essentially equivalent definitions of the automorphism. For example, in
paper~\cite{EvdP12}
it was defined as an automorphism of the Cayley scheme associated with the S-ring. In the S-ring language this definition is
equivalent to the following one.

\dfntn
A bijection $f:G\to G$ is an automorphism of an S-ring $\A$ if given $x,y\in G$ the basic sets containing
the elements $xy^{-1}$ and $x^f(y^f)^{-1}$ coincide. The group of all automorphisms of $\A$ is denoted by $\aut(\A)$.
\edfntn

It is easily seen that $G_{right}\le\aut(\A)$ and any basic set of $\A$ is invariant with respect to the group $\aut(\A)_1$.
Moreover, the group $\aut(\A)$ is the largest subgroup of $\sym(G)$ that satisfies these two properties.  One can see
that if $\A$ is a schurian S-ring  and $\Gamma\in\M(\A)$ where
$$
\M(\A)=\{\Gamma\le\sym(G):\ \Gamma\twoe\aut(\A)\ \,\text{and}\ \,G_{right}\le\Gamma\},
$$
then the S-ring associated with~$\Gamma$ is equal to~$\A$. It follows that  an S-ring $\A$ is schurian if and only if
$\S(\A)=\orb(\aut(\A)_1,G).$\medskip

Let $f\in\aut(\A)_1$.
Then any $\A$-set (in particular, $\A$-group)  is invariant with respect to the automorphism~$f$. 
Moreover,
for any $\A$-section $S$ we have $f^S\in\aut(\A_S)$. In particular, the S-ring $\A_S$ is schurian whenever so is~$\A$.\medskip 


The following result proved in \cite[Corollary~5.7]{EvdP12}  gives a criterion for the schurity of generalized wreath products that will be
 repeatedly used throughout the paper.

\thrml{060813a}
Let $\A$ be an S-ring over an abelian group~$G$. Suppose that $\A$ is an $S$-wreath product where~$S=U/L$. Then $\A$ is schurian if and only if
so are the S-rings $\A_{G/L}$ and $\A_U$ and there exist groups $\Delta_0\in\M(\A_{G/L})$ and $\Delta_1\in\M(\A_U)$ such that $(\Delta_0)^S=(\Delta_1)^S$.
\ethrm

The S-ring $\A$ is called {\it normal}, if $G_{right}$ is a normal subgroup of $\aut(\A)$, or, equivalently, if $\aut(\A)\le\hol(G)$. 
One can see that 
$\A$ is normal whenever it contains a normal S-ring over $G$ or is the tensor product of normal S-rings.
It was also proved in  \cite{EP01ce} that
when the group $G$  is cyclic, any normal S-ring over it is cyclotomic; in particular, any subgroup of $G$ is an $\A$-group.


\section{A sufficient condition of non-schurity}

Below given  a positive integer $n,$ we set
$$ \Omega^*(n) =
\css \Omega(n) & \mbox{if $n$ is odd}, \\  \Omega(n/2) & \mbox{if $n$ is even}, \ecss $$
where $\Omega(n)$ denotes the total number of prime divisors of $n$. It immediately follows from
the definition  that $\Omega^*(n) \le 1$  if and only if $n$ is a divisor of twice a prime number.
The goal of this section is the following sufficient condition of non-schurity.

\thrml{040211b}
Let $G_i$ be an abelian group of order $n_i,$ $\Omega^*(n_i)\ge 2$, $i=1,2$. Then
$G:=G_1\times G_2$ is not a Schur group.
\ethrm

\proof Suppose on the contrary that $G$ is a Schur group. Since the class of Schur groups is
closed with respect to taking subgroups, without loss of generality we can assume that
there exist odd primes $p_i$ and $q_i$ ($i=1,2$) so that
\qtnl{250512a}
n_i=p_iq_i\quad\text{or}\quad
n_i=4q_i\quad\text{or}\quad
n_i=8.
\eqtn
Thus we can assume that there exist groups $A_i$ and $B_i$ ($i=1,2$),  such that $1<A_i\le B_i<G_i$ and
\qtnl{200813a}
|G_i/A_i|\ge 3,\quad |B_i|\ge 3,\quad  |B_i/A_i|\le 2.
\eqtn
Indeed, if $n_i=p_iq_i$ or $n_i=4q_i$, set  $A_i$ to be the subgroup of $G_i$ of index $q_i$  and $B_i=A_i$, whereas if
$n_i=8$, set $A_i$ to be a  subgroup of $G_i$ of order $2$ and choose $B_i$ so that $|B_i:A_i|=2$.
A part of the subgroup lattice of~$G$ is given
in Figure~\ref{f1}.
\begin{figure}[h]
$\hspace{25mm}\xymatrix@R=10pt@C=20pt@M=0pt@L=5pt{
&  &  & \VRT{G} \ar@{-}[dl]\ar@{-}[dr]                     &  & &\\
  &  & \VRT{}\ar@{-}[dl]\ar@{-}[dr]                     &  &  \VRT{}\ar@{-}[dl]\ar@{-}[dr] & &\\
  & \VRT{}\ar@{-}[dl]\ar@{-}[dr] & & \VRT{}\ar@{-}[dl]\ar@{-}[dr] & & \VRT{}\ar@{-}[dl]\ar@{-}[dr] &\\
\VRT{G_1}\ar@{-}[dr] & & \VRT{}\ar@{-}[dl]\ar@{-}[dr] & & \VRT{}\ar@{-}[dl]\ar@{-}[dr] & & \VRT{G_2}\ar@{-}[dl]\\
& \VRT{B_1}\ar@{-}_{Q_1}[dr] & & \VRT{}\ar@{-}[dl]\ar@{-}[dr] & & \VRT{B_2}\ar@{-}^{Q_2}[dl] &\\
& &  \VRT{A_1}\ar@{-}[dr] & & \VRT{A_2}\ar@{-}[dl] & &\\
& & & \VRT{1} & & &\\
}$
\caption{}\label{f1}
\end{figure}

\medskip  We observe that by~\eqref{250512a} and~\eqref{200813a} each of the groups 
$$
H_1:=B_1,\quad H_2:=G_1/A_1,\quad H_3:=B_2,\quad H_4:=G_2/A_2,
$$
 is isomorphic either to $\mZ_r$ where $r$ is an odd prime or $r=4$, or to~$E_4$.
Denote by $\sigma_i$ the automorphism of $H_i$ that takes $x$ to~$x^{-1}$
in the former case, and by the involution $(x_1,x_2)\mapsto(x_2,x_1)$ in the latter case
(here $\sigma_i$ depends on the choice of generators). Then
$$
K_i:=\lg\sigma_i\rg\qaq K_{i,j}:=\lg \sigma_i\times\sigma_j\rg
$$
are subgroups of order~$2$ in the groups $\aut(H_i)$ and $\aut(H_i\times H_j)$ respectively.\medskip

Let us consider the cyclotomic S-ring $\cyc(K_i,H_i)$. One can see that the automorphism group of it  is a subgroup of the group $\hol(H_i)$:
this is obvious when $H_i=E_4$, whereas if $H_i$ is cyclic, then this follows because  it 
coincides with the automorphism group of the undirected cycle on~$H_i$ defined by~$K_i$.
Thus the S-ring $\cyc(K_i,H_i)$ is normal for all~$i$. Moreover, without loss of generality we can assume that
if $H_i=E_4$, then the unique proper $\A$-group in the S-ring $\cyc(K_i,H_i)$ coincides
with $A_1$ for $i=1$, with $B_1/A_1$ for $i=2$, with $A_2$ for $i=3$,  and with $B_2/A_2$
for $i=4$.\footnote{This can be done via an appropriate
decomposition $E_4=\mZ_2\times\mZ_2$.}\medskip

Set
\qtnl{240311z}
\A_1=\cyc(K_1\times K_3,H_1\times H_3),\qquad
\A_2=\cyc(K_{2,3},H_2\times H_3),
\eqtn
\qtnl{240311y}
\A_3=\cyc(K_{1,4},H_1\times H_4),\qquad
\A_4=\cyc(K_{2,4},H_2\times H_4).
\eqtn
Then each of these S-rings contains the tensor product of two normal S-rings. This implies
that $\A_i$ is normal for all~$i$. Set $Q_1=B_1/A_1$ and $Q_2=B_2/A_2$.
Then since $|Q_1|\le 2$ (see~\eqref{200813a}),  we conclude by statement~(2) of Lemma~\ref{050813b} that
$$
(\A_1)_{Q_1\times H_3}=\mZ Q_1\otimes\cyc(K_3,H_3)=(\A_2)_{Q_1\times H_3}
$$
and
$$
(\A_3)_{Q_1\times H_4}=\mZ Q_1\otimes\cyc(K_4,H_4)=(\A_4)_{Q_1\times H_4}.
$$
Thus one can form the following S-rings over $G_1\times H_3$ and $G_1\times H_4$:
$$
\A_{1,2}=\A_1\wr_{Q_1\times H_3}\A_2\qquad\text{and}\qquad \A_{3,4}=\A_3\wr_{Q_1\times H_4}\A_4.
$$
From the definition of the generalized wreath product it follows that $H_1$, $H_2$, $Q_2$, $H_3$
are  $\A_{1,2}$-sections, and $H_1$, $H_2$, $Q_2$, $H_4$ are $\A_{3,4}$-sections. Moreover, it is
easily seen that
$$
(\A_{1,2})_{G_1}= \cyc(K_1,H_1)\wr_{Q_1}\cyc(K_2,H_2)=(\A_{3,4})_{G_1}.
$$
Since $|Q_2|\le 2$  (see~\eqref{200813a}), we conclude by statement~(2) of Lemma~\ref{050813b} that
$$
(\A_{1,2})_S=(\cyc(K_1,H_1)\wr_{Q_1}\cyc(K_2,H_2))\otimes\mZ Q_2=(\A_{3,4})_S
$$
where $S=G_1\times Q_2$. Thus one can form an S-ring over the group $G$ as follows
$$
\A=\A_{1,2}\wr_S\A_{3,4}.
$$
To complete the proof it suffices to verify that this S-ring is not schurian.\medskip

Suppose on the contrary that $\A$ is schurian. Then by Theorem~\ref{060813a} the S-rings $\A_{1,2}$ and
$\A_{3,4}$ are schurian and  there exist groups $\Delta_{1,2}\in\M(\A_{1,2})$ and $\Delta_{3,4}\in\M(\A_{3,4})$
such that
$$
(\Delta_{1,2})^S=(\Delta_{3,4})^S.
$$
In particular, for any permutation $f_1\in\Delta_{1,2}$ leaving the
point~$1_{G_1\times H_3}$ fixed there exists a permutation $f_2\in\Delta_{3,4}$ leaving the
point~$1_{G_1\times H_4}$ fixed such that $f_1^S=f_2^S$. We claim that for any such $f_1$ we have
\qtnl{260512a}
(f_1)^{H_1\times H_3}\in\lg\sigma_1\times\sigma_3\rg.
\eqtn
However, if this is true, then the stabilizer of the point $1_{H_1\times H_3}$ in the group $(\Delta_{1,2})^{H_1\times H_3}$
is of order at most~$2$. Therefore any basic set of the S-ring associated with the latter group
is of cardinality $\le 2$. On the other hand, this S-ring coincides with $\A_1$
by the schurity of the S-ring $\A_{1,2}$ and the fact that the groups $\Delta_{1,2}$ and
$\aut(\A_{1,2})$ are 2-equivalent. However, by definition the S-ring~$\A_1$ has a basic set of
cardinality~$4$. Contradiction.\medskip

 To prove the claim let $f_{1,1}$ and $f_{1,2}$ be the automorphisms of the S-rings $\A_1$ and $\A_2$
induced by~$f_1$, and let $f_{2,3}$ and $f_{2,4}$ be the automorphisms of the S-rings $\A_3$ and $\A_4$
induced by~$f_2$. Then since these S-rings are normal, we have $f_{1,1}\in K_1\times K_3$,
$f_{1,2}\in K_{2,3}$, and $f_{2,3}\in K_{1,4}$, $f_{2,4}\in K_{2,4}$.
Clearly,
\qtnl{220311z}
(f_{1,1})^{Q_1\times H_3}=(f_{1,2})^{Q_1\times H_3}\quad\text{and}\quad
(f_{2,3})^{Q_1\times H_4}=(f_{2,4})^{Q_1\times H_4}
\eqtn
and due to the equality $(f_1)^S=(f_2)^S$ also
\qtnl{220311y}
(f_{1,1})^{H_1\times Q_2}=(f_{2,3})^{H_1\times Q_2}\quad\text{and}\quad
(f_{1,2})^{H_2\times Q_2}=(f_{2,4})^{H_2\times Q_2}.
\eqtn
Next, the permutations $f_{1,2}$, $f_{2,3}$ and $f_{2,4}$ are powers of the involutions
$\sigma_2\times\sigma_3$, $\sigma_1\times\sigma_4$ and $\sigma_2\times\sigma_4$ respectively.
Denote the corresponding exponents by $\varepsilon_{1,2},\varepsilon_{2,3}$ and $\varepsilon_{2,4}$
(all of them belong to  \{0,1\}).
Therefore, by the second equalities in~\eqref{220311z} and~\eqref{220311y} we have
\qtnl{0890812a}
\varepsilon_{1,2}=\varepsilon_{2,4}=\varepsilon_{2,3}.
\eqtn
Denote this number by $\varepsilon$. Then by the first equality in \eqref{220311y} the permutation
$(f_{1,1})^{H_1}$ is identical if and only if so is the permutation~$(f_{2,3})^{H_1}$, or equivalently
when $\varepsilon_{2,3}=1$. Similarly, by the first equality in~\eqref{220311z} the permutation
$(f_{1,1})^{H_3}$ is identical if and only if $\varepsilon_{1,2}=1$. Thus due to~\eqref{0890812a}
we have
$$
(f_1)^{H_1\times H_3}=f_{1,1}=(f_{1,1})^{H_1}\times(f_{1,1})^{H_3}
=\sigma_1^\varepsilon\times\sigma_3^\varepsilon,
$$
which proves~\eqref{260512a}. \bull

\section{Proofs  of Theorems~\ref {290713a} and~\ref{040213a}}
We begin with a simple observation that essentially restricts the structure of a Schur $p$-group.

\lmml{010813a}
Let $G$ be an abelian non-cyclic Schur $p$-group. Then $p<5$, and either $G=E_{p^k}$ or
$G=\mZ_{p^{}}\times\mZ_{p^k}$ for some~$k$.
\elmm
\proof Since $G$ is a non-cyclic Schur group, it follows from the P\"oschel theorem (see Introduction) that $p=2$ or
$3$. Denote by $r$ the rank of $G$. Then $G$ is elementary abelian whenever $r\ge 3$, for otherwise $G$ contains a subgroup isomorphic
to $\mZ_{p^2} \times E_{p^2}$ which is not a  Schur group  for $p=2$ by the CSS  and
for $p=3$ by Theorem~\ref{040211b}.  Let $r=2$. Then either $G$ is isomorphic to a group $\mZ_{p^{}}\times\mZ_{p^k}$ or contains a subgroup isomorphic to
$\mZ_{p^2}\times\mZ_{p^2}$. Since the latter is not a Schur group  (again, for $p=2$ by the CSS, and for $p=3$ by Theorem~\ref{040211b}), we
are done.\bull \medskip

\noindent{\bf Proof of Theorem~\ref{290713a}}. Let $G$ be an elementary abelian group of order $n=p^k$ where $p$ is a prime and $k\ge 2$.
Then the "if" part immediately follows from the CSS. To prove the "only if" part suppose that $G$ is a Schur group. Then
$p=2$ or $3$ by Lemma~\ref{010813a}.  Thus the required statement follows from Theorem~\ref{040211b} which shows that the groups $E_{2^k}$ and $E_{3^k}$
are not Schur  for $k\ge 6$ and  $k\ge 4$ respectively.\bull\medskip

\noindent{\bf Proof of Theorem~\ref{040213a}}.
Let $G$ be an abelian Schur group which is neither cyclic  nor elementary abelian. If $G$ is a $p$-group, then $G=\mZ_2\times\mZ_{2^k}$
or $G=\mZ_3\times\mZ_{3^k}$ by Lemma~\ref{010813a}, i.e. $G$ belongs to the first or sixth family. Thus in what follows we can assume that
$G$ is not a $p$-group. \medskip

Suppose first that there exist a group $H$, a prime $p\ge 3$ and an integer $k\ge 2$ such that
\qtnl{310512c}
G=H\times\mZ_{p^k}.
\eqtn
Then by Theorem~\ref{040211b} we have $\Omega^*(m) \le 1$ where $m=|H|$.
Thus $m$ divides $2 q$ for a prime~$q$. Since $G$ is not cyclic and $p \ge 3$, we conclude that
$m \in \{q,2q\}$ and either $q=2$ or $q=p$. It follows that $m\in\{2,4,p,2p\}$. Furthermore,
$H$ is not isomorphic to $\mZ_2$ or $\mZ_4$ because $G$ is not cyclic,
$H$ is not isomorphic to $\mZ_p$ because $G$ is not a $p$-group, and
$H$ is isomorphic to $\mZ_{2p}$ only if $p=3$ because $G$ has no
section  $E_{p^2}$ for $p \ge 5$ by Theorem~\ref{290713a}. Thus
\qtnl{310512b}
G=E_4\times\mZ_{p^k}\quad\text{or}\quad G=\mZ_{6}\times\mZ_{3^k}.
\eqtn
In particular, $G$ belongs to the third or seventh family. Thus, in what follows, without loss of generality we can assume that
$G$ is not of the form~\eqref{310512c}.\medskip

Let us turn to the general case. Since any subgroup of a Schur group is also Schur, the P\"oschel theorem implies
that the Sylow $p$-subgroup $G_p$ of $G$ is cyclic
for all prime $p \ge 5$. Since $G$ is not cyclic, this implies that so is not one of the groups $G_2$ and $G_3$.
Therefore,  either $E_4\le G$ or $E_9\le G$. Besides, by the CSS the group $E_4\times E_9$ is not Schur.
Thus to complete the proof it suffices to consider the following two cases.\medskip

\noindent {\bf Case 1:} $E_4\le G$ and $E_9\not\le G$.
The latter implies that the group $G_3$ is cyclic. Since $G_p$ is cyclic for $p \ge 5,$ 
the largest  group $G_{2'}\le G$ of odd order  is cyclic too.
Since also $G$ is not of the form~\eqref{310512c}, the number $|G_{2'}|$ is squarefree.
By Theorem~\ref{040211b} this implies that
$$
G = G_2 \times\mZ_{p^\varepsilon q^\delta}
$$
for some odd primes $p,q$ and $\varepsilon,\delta\in\{0,1\}$. Moreover,
$G_2$ is a Schur $2$-group. So $G_2 \le E_{32}$ or $G_2 \le\mZ_2\times\mZ_{2^k}$ by Lemma~\ref{010813a} and Theorem~\ref{290713a}.
Furthermore, $(\varepsilon,\delta)\ne (0,0)$ because $G$ is not a $2$-group, and also from  Theorem~\ref{040211b} it follows that $|G_2| \le 4$ when 
$(\varepsilon,\delta)=(1,1)$  and $G_2\le E_{16}$ when $(\varepsilon,\delta)=(0,1)$ or $(1,0)$. Thus  in any case
$$ 
G = E_4\times\mZ_{pq}\quad\text{or}\quad G\le E_{16}\times\mZ_p\quad\text{or}\quad
G\le\mZ_2\times\mZ_{2^k}\times\mZ_p.
$$ 
Thus $G$ belongs to the fourth family in the first case, to the third or fifth family in the second case and to
the second family in the third case.\medskip

\noindent {\bf Case 2:} $E_4\not\le G$ and $E_9\le G$.  
Reasoning as  in Case~1 we obtain that the group $G_{3'}$ is cyclic, where
$G_{3'}$ is the largest  subgroup of $G$ of order coprime to $3$.
Moreover, by Theorem~\ref{040211b} the number $|G_{3'}|$ divides $2 p$ for a prime $p \ne 3$. Thus
$$
G=G_3\times\mZ_{2^\varepsilon p^\delta}
$$
where $\varepsilon,\delta\in\{0,1\}$. On the other hand, $G_3$ is  a Schur $3$-group. Therefore
$G_3 \le E_{27}$ or $G_3=\mZ_3\times\mZ_{3^k}$ by  Lemma~\ref{010813a} and Theorem~\ref{290713a}.
Since $G$ is not of the form~\eqref{310512c},  in the latter case we have $k = 1$. Thus in any case $G_3 = E_9$ or $G_3 = E_{27}$.
Besides, by Theorem~\ref{040211b} in the latter case $(\varepsilon,\delta)=(1,0)$, and hence $G = E_9\times\mZ_6$
(the ninth family for $q=3$). Finally, if $G_3=E_9$, then
$$ 
G=E_9\times\mZ_2\quad\text{or}\quad
G=E_9\times\mZ_p\quad\text{or}\quad
G=E_9\times\mZ_{2p}.
$$ 
Thus $G$ belongs to the eighth family in the first or second case and to
the ninth family in the third case. Theorem~\ref{040213a} is completely proved. \bull

\section{Proof of Theorem~\ref{040213c}}

\sbsn\label{260813a}
Throughout this preliminary subsection $\A$ denotes an S-ring over a group
$G=H\times P$ where $H$ is an abelian group and $P=\mZ_p$ with prime~$p$
coprime to the order of~$H$.
In what follows given $X\subset G$ we set  $H_X=H\cap X$ and $H'_X=\pr_H(X)\setminus H_X$.

\lmml{010612a}
Given  $X\in\S(\A)$ the set $H'_X$ belongs to $\S^*(\A)$. Moreover, if $X$ meets both  $H$
and $G\setminus H$, then
\qtnl{010612b}
X=H_XP\cup H'_XP^\#\qaq X\cup H'_X=\pr_H(X)P.
\eqtn
\elmm
\proof Without loss of generality we can assume that $X\not\subset H$ (otherwise $H'_X=\emptyset$, $X$ does not meet $G\setminus H$, and both
statements are trivial). Set
$$
Y=(X^{[p]})^\sigma
$$
where $\sigma:x\mapsto x^m$ is an automorphism of $G$ such that $mp\equiv 1\pmod{n/p}$ with $n=|G|$. Then $Y$ consists of
all elements $x\in\pr_H(X)$ for which $xP\not\subset X$. On the other hand, from Theorems~\ref{261009b}
and \ref{261009w} it follows that $Y\in\S^*(\A)$. So either $X\subset Y$, or $X\cap Y=\emptyset$. Since also $Y\subset H$, our assumption
implies that $X\cap Y=\emptyset$, and so $H_X\cap Y=\emptyset$. It follows that
\qtnl{090813a}
H_X=\{x\in H:\ xP\subset X\}.
\eqtn
Thus $H'_X=Y$ which proves the first statement. To prove the second one suppose, in addition, that $X$ meets~$H$.
Since $H_X\ne\emptyset$, we conclude by Theorem~\ref{261009b} that
$$
X^{1\times\aut(P)}=X.
$$
Besides, by the transitivity of  $\aut(P)$ on $P^\#$ we have $x^{1\times\aut(P)}=xP^\#$ for all $x\in G$ with $x_P\ne 1$.
Therefore,
\qtnl{090813b}
H'_X=\{x\in H\setminus H_X:\ xP^\#\subset X\}.
\eqtn
 Thus the first equality in~\eqref{010612b} follows from~\eqref{090813a} and~\eqref{090813b}. The second equality immediately follows from the first one.\bull

\lmml{310512q}
Let  $H_1$ be the maximal $\A$-group contained in $H$. Suppose that $H_1\ne H$.
Then one of the following statements holds:
\nmrt
\tm{1} $\A=\A_{H_1}\wr\A_{G/H_1}$ with $\rk(\A_{G/H_1})=2$,
\tm{2} $\A$ is a $U/L$-wreath product  where $P \le L < G$ and $U=H_1L$.
\enmrt
\elmm
\proof Let $X$ be a basic set of $\A$ that meets $H\setminus H_1$. Then $X\not\subset H$, because otherwise $\lg X,H_1\rg$ is an $\A$-group
properly containing $H_1$, in contrast to the maximality of~$H_1$. Therefore $X$ meets both  $H$ and $G\setminus H$. By
Lemma~\ref{010612a} this implies that $H'_X\in\S^*(\A)$, and hence $Q:=\pr_H(X)P$ is an $\A$-set by the second equality in~\eqref{010612b}.
Suppose first that $Q=G$. Then the same equality implies that $X\supset G\setminus H$. Thus any basic set
that meets $H\setminus H_1$, contains $G\setminus H$. It follows that $X=G\setminus H_1$, i.e. statement~(1) holds.\medskip

Let now  $Q<G$. Then $\rad(Q)$ is an $\A$-group such that  $P\le \rad(Q)<G$. Therefore the minimal $\A$-group
$L$ that contains~$P$, does not equal~$G$. To complete the proof it suffices to verify that
\qtnl{040213d}
P\le \rad(Y) \,\quad Y\in\S(\A)_{G\setminus U}
\eqtn
where $U=H_1L$. Indeed, in this case the minimality of $L$ implies that $L\le \rad(Y)$ for all $Y$,  and hence
the S-ring $\A$ is the $U/L$-wreath product as required.\medskip

To prove~\eqref{040213d} it suffices to verify that any basic set $Y\subset G\setminus U$ that meets $H$
is of the form
\qtnl{120713l}
Y=\pr_H(Y)P
\eqtn
(indeed, in this case obviously $P\le\rad(Y)$ and the union of all these $Y$'s coincides with $G\setminus U$).
To verify the claim let $Y$ be as above. Then $Y$ does not meet $H_1\le U$, and so $Y\not\subset H$ by the maximality of $H_1$.
It follows that $Y$ meets both  $H$ and $G\setminus H$. By  Lemma~\ref{010612a} applied to $X=Y$, this implies that $H'_Y$
is an $\A$-set, and hence $H'_Y\subset H_1$ by the maximality of $H_1$, whence $H'_YP  \subset H_1P \subset U$.
It follows that $\pr_H(Y)=H_Y$, for otherwise the set $H'_YP\cap Y$ is not empty, and so $Y$ does meet $U$, which is impossible
by the choice of~$Y$. This implies that the set $H'_Y$ is empty. Thus  \eqref{120713l} follows from~\eqref{010612b} 
with $X$ replaced by~$Y$.\bull

\sbsn Let us turn to the proof of the theorem.  Keeping the notation of Subsection~\ref{260813a} we assume that
$H=E_4$ and  $p$ is an odd prime (if $p=2$, then $G=E_8$ is a Schur group by Theorem~\ref{290713a}).
Then any proper section of $G$ is isomorphic to the group $E_4$ or a cyclic group of order $2$, $p$ or $2p$. Since all
these groups are Schur (see Introduction), we conclude that the S-ring $\A_S$ is schurian for any proper $\A$-section~$S$. 
The rest of the proof is divided into three mutually exclusive cases. \medskip

\noindent{\bf Case 1:} $H$ is not an $\A$-group.  Then the hypothesis of Lemma~\ref{310512q} is satisfied.
So one of two statements of that lemma holds. If statement~(1) holds, then 
$H_1=1$ or $H_1=\mZ_2$. In the former case $\A$ is schurian because $\rk(\A)=2$. In the latter case
the S-ring $\A$  is schurian as  the wreath product of schurian S-rings $\A_{H_1}$ and $\A_{G/H_1}$.\medskip

To complete the case suppose that statement (2) of Lemma~\ref{310512q} holds. Then $\A$ is the $U/L$-wreath product
where $P \le L < G$ and $U=H_1L$. Suppose first that the wreath product is not proper. Then $U=G$, and hence
$H_1=\mZ_2$ and $G=H_1\times L$. By statement~(2) of Lemma~\ref{050813b} this implies that $\A=\A_{H_1}\otimes\A_L$.
Thus the S-ring $\A$  is schurian as  the tensor product of schurian S-rings.\medskip

In the remaining case $\A$ is a proper $U/L$-wreath product.  Then $|U/L|\le 2$ because $|L|\ge p$ and $|G/U|\ge 2$.
Furthermore, since the S-rings $\A_U$ and $\A_{G/L}$ are schurian,
the groups  $\Delta_0=\aut(\A_{G/L})$ and  $\Delta_1=\aut(\A_U)$ belong to the sets $\M(\A_{G/L})$ and $\M(\A_U)$
respectively. Finally,
$$
(\Delta_1)^{U/L}=(\Delta_0)^{U/L},
$$
because each of these two groups is of order  $|U/L|\le 2$.
Thus the S-ring  $\A$ is schurian by Theorem~\ref{060813a}.\medskip

\noindent{\bf Case 2:} $H$ is an $\A$-group while $P$ is not.\footnote{In fact, the case when $P$ is not an $\A$-group is
dual to Case 1 in the sense of the duality theory of S-rings over an abelian group, see~\cite{EvdP09}; here we prefer
a direct proof.}  It suffices to verify that there exists a non-identity $\A$-group $L \le H$ such that
\qtnl{050213b}
\A=\A_H \wr_{H/L} \A_{G/L}.
\eqtn
Indeed, in this case the fact that the S-ring $\A$ is schurian can be proved in the same way as in the end of Case~1.\medskip

To prove \eqref{050213b} we observe that given $X\in\S(\A)$ there exists a positive integer $m=m(X)$ such that
$|X\cap C|=m$ for any $H$-coset $C$ intersecting $X$ (Lemma~\ref{090608a}). We claim that
\qtnl{050713a}
m(X)\in\{2,4\},\quad X\in\S(\A)_{G\setminus H}.
\eqtn
Indeed, suppose on the contrary that the number $m(X)$ is odd for some $X\in\S(\A)_{G\setminus H}$. Then  $Y=X^{[2]}$ is
a nonempty subset of $P^\#$ where $X^{[2]}$ is the set defined in~\eqref{030713a}. By Theorem~\ref{261009w}
we have $Y\in\S^*(\A)$. Since $|P|$ is prime, it follows that $P=\lg Y\rg$ is an $\A$-group. Contradiction.\medskip

To complete the proof of~\eqref{050213b} it suffices to find a non-identity group $L \le H$ such that $\rad(X)=L$ for
all $X\in\S(\A)_{G\setminus H}$. However, any nonempty subset of  $H$ of even cardinality is a coset by a subgroup of $H$,
and if the subset is of cardinality $2$, then its complement in $H$ is a coset by the same subgroup. Therefore
due to \eqref{050713a}, it suffices to verify that
\qtnl{050713b}
\pr_H(X\cap C)=\pr_H(Y\cap D)\quad\text{or}\quad \pr_H(X\cap C)\,\cap\,\pr_H(Y\cap D)=\emptyset
\eqtn
for all $X,Y\in\S(\A)_{G\setminus H}$ and all $C,D\in (G/H)^\#$. To prove \eqref{050713b}, suppose that the sets  $\pr_H(X\cap C)$ and
$\pr_H(Y\cap D)$ intersect. Then there exists an element $\sigma\in 1\times\aut(P)$ such that $X^\sigma\cap Y\ne\emptyset$
(here we used the fact that $\aut(P)$ acts transitively on $P^\#$). By Theorem~\ref{261009b} the set $X^\sigma$ belongs
to $\S(\A)$. Thus $X^\sigma=Y$and $C^\sigma=D$, and hence
$$
\pr_H(X\cap C)=\pr_H(X^\sigma\cap C^\sigma)=\pr_H(Y\cap D)
$$
as required.\medskip

\noindent{\bf Case 3:} $H$ and $P$ are $\A$-groups. Then the S-rings $\A_H$ and $\A_P$ are cyclotomic: the former is 
easy~\footnote{Any S-ring over $H$ is of the form $\cyc(K,H)$ where $K=\aut(H)\cong\sym(3)$ or $K$ is a subgroup of order~$2$ in $\aut(H)$ or $K=1$.}
whereas the latter follows from~\cite{EP01ce}.
Set
$$
K=K_1\times K_2
$$
where $K_1$ and $K_2$ are the unique cyclic subgroups of $\aut(H)$ and $\aut(P)$  such that
\qtnl{110713j}
\A_H=\cyc(K_1,H)\qaq \A_P=\cyc(K_2,P).
\eqtn
Without loss of generality we can assume that $\A\ne\A_H\otimes\A_P$, for otherwise 
the S-ring $\A$ is schurian as  the tensor product of schurian S-rings. Then $K_1\ne 1$ by statement~(2) of Lemma~\ref{050813b}. 
It follows that  $|K_1|\in\{2,3\}$, and the group $K_1$
has a unique regular orbit on $H$; denote it by~$A$. Let us call a basic set $X\in\S(\A)$ {\it highest} if $X_H=A$ and $X_P\ne 1$ where
$X_H=\pr_H(X)$ and $X_P=\pr_P(X)$. By statement~(1) of Lemma~\ref{050813b} highest basic sets do exist.

\lmml{100713b}
In the above notations the following statements hold:
\nmrt
\tm{1} for any highest basic set $X$ we have $X_P\in\orb(K_2,P)$,
\tm{2} the group $1\times\aut(P)$ acts transitively on the highest basic sets,
\tm{3} if $X$ is a non-highest basic set, then $X=X_H\times X_P$, and one of the sets $X_H$ or $X_P$ is a singleton.
\enmrt
\elmm
\proof Statements (1) and (3) are obvious. Statement~(2) follows from Theorem~\ref{261009b} because the group $1\times\aut(P)$ is contained
in the center of~$\aut(G)$ and $\aut(P)$ acts transitively on $P^\#$.\bull\medskip

Let $X$ be a highest basic set of $\A$. Then $|X|$ divides $|A||K_2|$ by statement~(2) of Lemma~\ref{100713b}. On the other hand,
from Lemma~\ref{090608a} it follows that $|X_P|$ divides $|X|$. Since also $|X_P|=|K_2|$, we conclude that
$|K_2|$ divides $|X|$. Since $|A|\in\{2,3\}$, we conclude that
$$
|X|=|K_2|\quad\text{or}\quad |X|=|A||K_2|.
$$
However, in the latter case  by Lemma~\ref{100713b}
any basic set of $\A$ is the product of its projections, which is impossible because $\A\ne\A_H\otimes\A_P$
(see above). Thus
\qtnl{110713t}
|X|=|K_2|.
\eqtn
This implies that in the decomposition
\qtnl{110713s}
X=\bigcup_{a\in A}\{a\}\times X_a
\eqtn
where $X_a=\{x\in P:\ (a,x)\in X\}$, the sets $X_a$ are pairwise disjoint. Moreover, from Lemma~\ref{090608a} applied to $H=P$ it follows that
these sets have the same size. Taking into account that their union is equal to $X_P$, we conclude  by~\eqref{110713t} that
\qtnl{280813a}
|X_a|=|K_2|/|A|,\quad a\in A.
\eqtn
On the other hand, from statement~(2) of  Lemma~\ref{100713b} and the fact that the sets in the right-hand side of~~\eqref{110713s} are disjoint it follows
that  for each $a\in A$ the setwise stabilizer of $\{a\}\times X_a$  in the group $1\times K_2$  is of order $|X_a|$. By~\eqref{280813a}
this implies that the stabilizer of $X_a$ in $K_2$ is equal to the subgroup $M$ of  index $|A|$ in $K_2$. Thus
\qtnl{110713a}
\{X_a:\ a\in A\}=\orb(M,X_P)
\eqtn
for all highest basic sets $X$ of~$\A$.\medskip

To complete the proof we will show that $\A=\cyc(K_0,G)$ where
$$
K_0=\{(\sigma_1,\sigma_2)\in K_1\times K_2:\ \pi_1(\sigma_1)=\pi_2(\sigma_2)\}
$$
with $\pi_2:K_2\to K_2/M$ the quotient epimorphism, and $\pi_1:K_1\to K_2/M$ an isomorphism defined as follows. Let $X$ be a highest basic set of $\A$.
Then from statement~(1) of Lemma~\ref{100713b} and equality~\eqref{110713a} it follows that $K_2$ acts on the set $\{X_a:\ a\in A\}$
as a regular cyclic group, and the kernel  of this action equals~$M$.  Since $|K_2/M|\in\{2,3\}$, it follows that given $\sigma_1\in K_1$ there exists
a unique $M$-coset $C$ in $K_2$  such that for any $\sigma_2\in C$ we have
$$
X_{a^{\sigma_1}}=(X_a)^{\sigma_2}=(X_a)^C.
$$
It immediately follows that $\pi_1:\sigma_1\mapsto C$ is the required isomorphism and $\pi_1(\sigma_1)=\pi_2(\sigma_2)$ for all $\sigma_2\in C$.\medskip

It easily follows from \eqref{110713s} and the definition of $K_0$
that the set $X$ is $K_0$-invariant. Moreover, since the group $K_0\le\aut(G)$ acts regularly on $X_H$ and $X_P$, 
the action of it on $X$ is faithful and fixed point free. So
$|X|\ge |K_0|$. Furthermore, from the definition of $K_0$ it follows that $|K_0|=|K_2|$. Thus due to~\eqref{110713t} we conclude that $|X|=|K_0|$, and hence
\qtnl{110713g}
X\in\orb(K_0,G).
\eqtn
If $Y$ is another  highest basic set of $\A$, then by statement~(2) of Lemma~\ref{100713b} we have $Y=X^\sigma$ for some $\sigma\in 1\times\aut(P)$. However,
the permutation $\sigma$ obviously centralizes the group $K_0$. Therefore $Y$ is also an orbit of $K_0$, and hence relation
\eqref{110713g} holds for all highest basic sets. It remains to note that due to~\eqref{110713j}, statement~(3) of Lemma~\ref{100713b} and the
definition of~$K_0$, the relation also holds for all non-highest basic sets.\bull

\end{document}